\documentclass[10pt]{amsart}

\usepackage{amsfonts,amssymb}

\newtheorem{thm}{Theorem}[section]
\newtheorem{prop}[thm]{Proposition}
\newtheorem{cor}[thm]{Corollary}

\newtheorem*{thmA}{Theorem A}
\newtheorem*{thmB}{Theorem B}

\theoremstyle{remark}

\theoremstyle{definition}

\newcommand*\isom{%
  \xrightarrow{\sim}%
}

\def\PP{\mathbb{P}}

\def\zz{\mathbb{Z}}
\def\cc{\mathbb{C}}

\def\HH{\mathbb{H}}
\def\deldelbar{\partial \bar{\partial}}

\numberwithin{equation}{section}

\begin{document}

\title[Asymptotic behavior of the Kawazumi-Zhang invariant]{Asymptotic behavior of the Kawazumi-Zhang invariant for degenerating Riemann surfaces}

\author{Robin de Jong}



\begin{abstract} Around 2008 N. Kawazumi and S. Zhang introduced a new fundamental  numerical invariant for compact Riemann surfaces. One way of viewing the Kawazumi-Zhang invariant is as a quotient of two natural hermitian metrics with the same first Chern form on the line bundle of holomorphic differentials. In this paper we determine precise formulas, up to and including constant terms, for the asymptotic behavior of the Kawazumi-Zhang invariant for degenerating Riemann surfaces. As a corollary we state precise asymptotic formulas for the beta-invariant introduced around 2000 by R. Hain and D. Reed. These formulas are a refinement of a result Hain and Reed prove in their paper. We illustrate our results with some explicit calculations on degenerating genus two surfaces. 
\end{abstract}

\maketitle

\thispagestyle{empty}

\section{Introduction}

Consider a compact and connected Riemann surface $M$ of genus $h \geq 1$. In this paper we will be concerned with the \emph{canonical} or \emph{Bergman}  K\"ahler form $\mu$ on $M$, given as follows: let $(\omega_1,\ldots,\omega_h)$ be an orthonormal basis of holomorphic differentials on $M$, endowed with the $L^2$-inner product given by
\begin{equation} \label{inner}
<\alpha,\beta> = \frac{\sqrt{-1}}{2} \int_M \alpha \, \bar{\beta} \, .
\end{equation}
Then we put $\mu = \frac{\sqrt{-1}}{2h} \sum_{i=1}^h \omega_i \, \bar{\omega}_i$. Note that we have $\int_M \mu =1$ so that we may view $\mu$ as a probability measure on $M$. The canonical K\"ahler form $\mu$ and the resulting spectral theory of $M$ are useful in a variety of applications, ranging from arithmetic geometry \cite{ar} \cite{ce} \cite{fa} \cite{la} \cite{so}, to string perturbation theory \cite{abmnv} \cite{amv} \cite{bm} \cite{bk} \cite{dp88}.

Let $\Delta_\mu$ be the Laplace operator associated to $\mu$, i.e. the endomorphism of $L^2(M,\mu)$ determined by the condition
\[ \Delta_\mu(f) \cdot \mu = \frac{\deldelbar}{\pi \sqrt{-1}} f \]
for all $C^\infty$ functions $f$ in $L^2(M,\mu)$. Let $\lambda_0<\lambda_1\leq \lambda_2 \leq \ldots$ be the eigenvalues of $\Delta_\mu$ and let $(\phi_\ell)_\ell$ for $\ell=0,1,\ldots$ be a corresponding orthonormal basis of real eigenfunctions in $L^2(M,\mu)$. The invariant mentioned in the title refers to the invariant
\begin{equation} \label{defphi} \varphi(M) = \sum_{\ell>0} \frac{2}{\lambda_\ell} \sum_{m,n=1}^h \left| \int_M \phi_\ell \, \omega_m  \bar{\omega}_n \right|^2
\end{equation}
associated to $M$. This fundamental invariant was introduced independently by N. Kawazumi \cite{kaw} and S. Zhang \cite{zh}, around 2008, from different perspectives. In \cite{kaw} it arises in the context of a study of the behavior of the Johnson homomorphism on the moduli space of Riemann surfaces. The motivation in \cite{zh} comes from arithmetic geometry, where $\varphi$ occurs as a local contribution in a formula relating the height of the so-called Gross-Schoen cycle on a curve over a number field to the self-intersection of its relative dualizing sheaf. One way of viewing $\varphi(M)$ is as a quotient of two natural hermitian metrics, with the same first Chern form, on the line bundle of holomorphic differentials on $M$ \cite{kaw}.

It is easy to verify from the definition (\ref{defphi}) that $\varphi(M)=0$ if $h=1$, and that $\varphi(M)>0$ if $h \geq 2$. Our goal in this paper is to study the asymptotic behavior of $\varphi$ in a family of Riemann surfaces degenerating into a stable curve with a single node. There are two cases to be dealt with: one where after removing the node a disconnected surface remains (the separating case), and one where after removing the node a connected surface remains (the non-separating case). We use in both cases a standard degenerating model over the punctured unit disk described by J. Fay in Chapter III of his book `Theta functions on Riemann surfaces' \cite{fay}. 

The behavior of spectral invariants on degenerating Riemann surfaces has mainly been studied in the context of the metric of constant curvature $-1$, see \cite{hej} \cite{ji1} \cite{ji2} \cite{wo1} \cite{wo2} for example. The present work, in the context of the Bergman metric, extends the picture laid down in the papers \cite{jo} and \cite{we} by J. Jorgenson and R. Wentworth on the canonical Green's function and the Faltings delta-invariant, and the paper \cite{hrar} by R. Hain and D. Reed on the beta-invariant.

Our main result is as follows.
\begin{thmA}
Let $M_t$ be a Fay's family of compact Riemann surfaces of genus $h$, degenerating as $t\to 0$ to the union of two surfaces $M_1,M_2$ of genera $h_1,h_2 \geq 1$, joined at a single node. Then the formula
\[ \lim_{t\to 0} \left[ \varphi(M_t) + \frac{2h_1h_2}{h} \log |\tau| \right] = \varphi(M_1) + \varphi(M_2) \]
holds.
For a Fay's family $M_t$ of genus $h+1$ degenerating to a surface $M$ of genus $h \geq 1$ with two distinct points $a,b$ identified at a node, the formula
\[ \lim_{t\to 0} \left[ \varphi(M_t) + \frac{h}{6(h+1)} \log |\tau| \right] = \varphi(M) - \frac{5h}{3(h+1)} g(a,b)  \]
holds. Here $g$ is the canonical Green's function on the surface $M$.
\end{thmA}
For the precise relation between the parameters $\tau$ and $t$ and for the notion of canonical Green's function we refer to the text below. 

Before we discuss applications of Theorem A we would like to point out the similarity of our result with the main result of \cite{we}. In this paper R. Wentworth studies the asymptotic behavior of the canonical Green's function and the Faltings $\delta$-invariant \cite{fa} for degenerating Riemann surfaces, up to and including constant terms. He obtains the following limit formulas:
\[ \lim_{t\to 0} \left[ \delta(M_t) + \frac{4h_1h_2}{h} \log |\tau| \right] = \delta(M_1) + \delta(M_2) \]
in the separating case, and
\begin{align*}
\lim_{t\to 0} & \left[ \delta(M_t) + \frac{4h+3}{3(h+1)} \log |\tau| + 6 \log(-\log|\tau|) \right] \\
&
= \delta(M) - \frac{2(2h-3)}{3(h+1)} g(a,b) - 2 \log (2\pi)
\end{align*}
in the non-separating case. Our derivation of Theorem A will be heavily based on the detailed results in \cite{we} leading to these limit formulas.

In \cite{zh} Zhang also studies the invariant
\begin{equation} \label{deflambda} \lambda(M) = \frac{h-1}{6(2h+1)}\varphi(M) + \frac{1}{12}\delta(M) - \frac{h}{3} \log(2\pi)
\end{equation}
of the Riemann surface $M$. Note that in \cite{zh}, the delta-invariant differs by an additive constant $4h\log(2\pi)$ from the one used here and in \cite{we}. In the paper \cite{djsecond} it is shown that $(8h+4)\lambda$ is a natural representative of the $\beta$-invariant defined, up to constants, on the moduli space of compact Riemann surfaces by Hain and Reed \cite{hrar}. The $\beta$-invariant arises naturally in the context of a study of normal functions associated to Ceresa cycles on jacobians. The main result of \cite{hrar} states asymptotic formulas for $\beta$, but only up to $\log|\tau|$ terms. 

Putting $\beta=(8h+4)\lambda$ and combining the asymptotics of $\delta$ and $\varphi$ stated above we obtain from (\ref{deflambda}) the following more precise result.
\begin{thmB}
Let $M_t$ be a Fay's family of compact Riemann surfaces of genus $h$, degenerating as $t\to 0$ to the union of two surfaces $M_1,M_2$ of genera $h_1,h_2 \geq 1$, joined at a single node. Then for Hain-Reed's $\beta$-invariant, the formula
\begin{align*} \lim_{t\to 0} & \left[ \beta(M_t) + 4h_1h_2 \log |\tau| \right] = \frac{2}{3}(h-1)(\varphi(M_1)+\varphi(M_2))\\ & +\frac{2h+1}{3}(\delta(M_1)+\delta(M_2)) - \frac{(8h+4)h}{3}\log(2\pi)
\end{align*}
holds.
For a Fay's family $M_t$ of genus $h+1$ degenerating to a surface $M$ of genus $h \geq 1$ with two distinct points $a,b$ identified at a node, the formula
\begin{align*}
\lim_{t\to 0} & \left[ \beta(M_t) + (h+1) \log |\tau| + 2(2h+3) \log(-\log|\tau|) \right] \\ & = \frac{2h}{3}\varphi(M) +\frac{2h+3}{3}\delta(M) - 2(h-1) g(a,b) -
\left(\frac{8h(h+3)}{3}+6\right)\log(2\pi)
\end{align*}
holds. Here $g$ is the canonical Green's function on the surface $M$.
\end{thmB}
Theorem B improves upon the main result of \cite{hrar} in that we now also have some control over the constant terms. This may have applications in a precise study of limits of periods and mixed Hodge structures associated to the Ceresa cycle. 

Another application of our main result is to determine precise constants in explicit formulas involving $\varphi$. An illustration of this is given in \cite{djspecial} where the `energy' of the Weierstrass points on a hyperelliptic Riemann surface $M$ is expressed in terms of $\varphi(M)$.

The contents of this paper are as follows. In Section \ref{prelim} we present a collection of formulas that we will use as a reference chart. In the next two sections we prove Theorem A, treating the separating case first in Section \ref{sep}, and then the non-separating case in Section \ref{nonsep}. We illustrate our results by investigating the case of degenerating genus two surfaces in more detail in Section \ref{example}. This features a precise asymptotic analysis of the well known Siegel modular form $\chi_{10}$.

\section{Preliminaries} \label{prelim}

As before let $M$ be a compact and connected Riemann surface of genus $h \geq 1$, and let $\mu$ be its canonical K\"ahler form. The canonical Green's function is the generalized function on $M \times M$ determined by the conditions
\begin{equation}
\partial_y \bar{\partial}_y \, g(x,y) =\pi \sqrt{-1} \, (\mu(y) - \delta_x)
\end{equation}
and
\begin{equation} \label{normalization}
\int_M g(x,y) \,\mu(y) = 0
\end{equation}
for all $x$ in $M$. One readily derives the symmetry property
\begin{equation} \label{symmetry}
g(x,y) = g(y,x)
\end{equation}
for all distinct $x,y$ in $M$, as well as the local expansion
\begin{equation}
g(x,y) = \log |z(x) - z(y)| + O(|z(x)-z(y)|)
\end{equation}
for all distinct $x,y$ in a coordinate chart $z \colon U \isom D$ of $M$. Here $D$ denotes the open unit disk in the complex plane. The $O$-term is a $C^\infty$ function.

The canonical Green's function $g$ induces a natural hermitian metric $\|\cdot \|$ on the holomorphic line bundle $L=\mathcal{O}_{M \times M}(\Delta)$ on $M \times M$, where $\Delta$ is the diagonal, by putting $\log \|1\|(x,y)=g(x,y)$ for distinct $x,y$ in $M$. Here $1$ denotes the canonical generating section of $\mathcal{O}_{M \times M}(\Delta)$. We write $k$ for the first Chern form of $(L,\|\cdot\|)$. Note that it represents the class of $\Delta$ in $\mathrm{H}_\mathrm{dR}^2(M \times M)$. The Green's function $g$ essentially inverts $\Delta_\mu$, in the sense that
\[ f(x) = -\int_M g(x,y) \, \Delta_\mu(f)(y) \, \mu(y) + \int_M f \, \mu \]
holds for all $x \in M$ and all $C^\infty$ functions $f$ on $M$. 

Let $(\omega_1,\ldots,\omega_h)$ be an orthonormal basis of holomorphic differentials on $M$, equipped with the $L^2$-inner product (\ref{inner}). Let $\pi_1,\pi_2$ denote the projections of $M \times M$ on the first and second factor, respectively. Then we have the explicit formula \cite{ar} \cite{zh}
\begin{equation} \label{kDeltaom}
k = \pi_1^* \mu + \pi_2^* \mu - \sqrt{-1} \sum_{i=1}^h (\omega_i(x) \bar{\omega}_i(y)+\omega_i(y) \bar{\omega}_i(x)) \,
\end{equation}
Write $\nu$ for the $(2,2)$-form $k^2$ on $M \times M$. A small calculation yields
\begin{equation} \label{kDeltasq}
\nu = 2 \, \pi_1^* \mu \, \pi_2^* \mu - \sum_{i,j=1}^h (\omega_i(x) \bar{\omega}_j(x) \bar{\omega}_i(y) \omega_j(y) + \bar{\omega}_i(x) \omega_j(x) \omega_i(y) \bar{\omega}_j(y)) \, .
\end{equation}
Note that the $(2,2)$-form $\nu$ integrates to the topological Euler characteristic of $M$:
\begin{equation} \label{area}
\int_{M \times M} 
\nu = \Delta . \Delta = 2-2h \, .
\end{equation}
From the explicit formula (\ref{kDeltasq}) we further deduce that for $p$ fixed
\begin{equation} \label{fiberintproj}
\int_{M \times M} g(x,p) \, \nu(x,y) = 0 \, . 
\end{equation}
For our analysis of $\varphi(M)$ we will make use of the identity
\begin{equation} \label{phikDelta}
\varphi(M) = \int_{M \times M} g(x,y) \, \nu (x,y) \, , 
\end{equation}
for a proof of which we refer to \cite{zh}, Proposition 2.5.3.

Let $(A_1,\ldots,A_h,B_1,\ldots,B_h)$ be a symplectic basis of $\mathrm{H}_1(M,\zz)$. 
Let $(v_1,\ldots,v_h)$ be a corresponding normalized basis of holomorphic differentials:
\[ \int_{A_i} v_j = \delta_{ij} \, , \quad \int_{B_i} v_j = \Omega_{ij} \, . \]
We call $\Omega=\Omega_{ij}$ the period matrix of $(v_1,\ldots,v_h)$ with respect to the chosen homology basis. From Riemann's bilinear relations we have that $\Omega_{ij}$ is symmetric, and $\mathrm{Im} \, \Omega_{ij}$ is positive definite. They also yield the expression \cite{we}, Section~2 
\begin{equation} \label{muv}
 \mu = \frac{ \sqrt{-1}}{2h} \sum_{i,j=1}^h (\mathrm{Im} \, \Omega)_{ij}^{-1} v_i \, \bar{v}_j \, .
\end{equation}
In a similar way one derives from (\ref{kDeltaom}) the expression
\begin{equation} \label{kDeltav} k = \pi_1^* \mu + \pi_2^* \mu - \sqrt{-1} \sum_{i,j=1}^h (\mathrm{Im} \, \Omega)_{ij}^{-1} (v_i(x) \bar{v}_j(y) + v_i(y)\bar{v}_j(x) )\, .
\end{equation}

\section{Separating case} \label{sep}

We begin by introducing the degeneration model for a family of Riemann surfaces developing a single separating node as described by J. Fay, \cite{fay} Chapter III. The basic references for this section are \cite{fay} and \cite{we}, esp. Sections 3 and 6 of the latter. 

Let $M_1,M_2$ be two compact Riemann surfaces of genus $h_1,h_2 \geq 1$, respectively, and choose points $p_1 \in M_1$ and $p_2 \in M_2$. Let $z_i \colon U_i \isom D$ for $i=1,2$ be local coordinate charts centered at the $p_i$, with $D$ denoting the open unit disk in the complex plane. Let $\dot{D}$ denote the punctured unit disk; then we put $S = \{ (x,y,t) \in \dot{D}^3 : xy=t \}$. For each $t \in \dot{D}$ we remove from $M_i$ the set $|z_i| \leq |t|$ for $i=1,2$, and we glue the annuli $\{ |t| < |z_i| <1 \}$ on the resulting open surfaces along the maps given by
\[ z_1 \mapsto (z_1, t/z_1,t) \in S \, , \quad z_2 \mapsto (t/z_2,z_2,t) \in S \, . \]
We obtain a family of compact Riemann surfaces $M_t$ over $\dot{D}$ of genus $h=h_1+h_2$. This family extends to a proper analytic family $M \to D$ whose central fiber $M_0$ is the union of $M_1,M_2$ with the points $p_i$ identified.

For each $t \in D$ we choose a symplectic basis of $\mathrm{H}_1(M_t,\zz)$ by extending across $(M_1 \setminus U_1) \times \dot{D} \subset M$ and $(M_2 \setminus U_2) \times \dot{D} \subset M$ canonical homology bases of $M_1,M_2$ lying in $M_1 \setminus U_1$ and $M_2 \setminus U_2$, respectively. In \cite{fay}, Proposition 3.1 we find the following key result.
\begin{thm} \label{faysep}
(Fay) For $t \neq 0$ sufficiently small, the elements of the normalized basis of holomorphic differentials $v_1(x,t),\ldots,v_h(x,t)$ on $M_t$ corresponding to the chosen homology basis have the following expansions with respect to the parameter $t$: for $1 \leq i \leq h_1$,
\[ v_i(x,t) = \left\{ \begin{array}{ll}
v_i^{(1)}(x) + t\, v_i^{(1)}(p) \omega^{(1)}(x,p) + o(t) & x \in M_1 \setminus U_1 \\
-t \, v_i^{(1)}(p) \omega^{(2)}(x,p) + o(t) & x \in M_2 \setminus U_2 \, ; \end{array} \right. \]
and for $h_1+1 \leq j \leq h$,
\[ v_j(x,t) = \left\{ \begin{array}{ll}
v_j^{(2)}(x) + t \, v_j^{(2)}(p) \omega^{(2)}(x,p) + o(t) & x \in M_2 \setminus U_2 \\
-t \, v_j^{(2)}(p) \omega^{(1)}(x,p) + o(t) & x \in M_1 \setminus U_1 \, . \end{array} \right. \]
Here $v_i^{(1)}(x)$ and $v_j^{(2)}(x)$ are normalized bases of holomorphic differentials on $M_1,M_2$ with respect to the chosen homology bases, and $\omega^{(i)}(x,y)$ are the canonical differentials of the second kind on $M_i \times M_i$ for $i=1,2$. The point $p$ denotes the point $p_i$ on either $M_i$, and the evaluations of the differentials are carried out in the local coordinates $z_i$. Each $o(t)$ term is a holomorphic differential on $M_t$, with the property that $\lim_{t \to 0} o(t)/t^2$ is a meromorphic differential on $M_1 \cup M_2$ with a single pole at $p$.
\end{thm}
\begin{cor} \label{periodsep} Let $\Omega_1,\Omega_2$ be the period matrices of the $v_i$ on $M_1$ and of the $v_j$ on $M_2$. Let $R$ be the vector
\[ R = (v_1(p_1),\ldots,v_{h_1}(p_1),v_{h_1+1}(p_2),\ldots,v_h(p_2)) \in \cc^h \, . \]
Then the period matrix $\Omega_t$ of the $v_i(x,t)$ has the expansion
\[ \Omega_t = \left( \begin{array}{cc}
\Omega_1 & 0 \\ 0 & \Omega_2 \end{array} \right) + 2\pi \sqrt{-1} t \, {}^t R \cdot  R + o(t) \]
as $t \to 0$ over $\dot{D}$.
\end{cor}
\begin{proof}
See \cite{fay}, Corollary 3.2.
\end{proof}
Theorem \ref{faysep} and Corollary \ref{periodsep} lie at the basis of the proofs of the propositions below. For each $i=1,2$ write $M_{t,i}=M_i \setminus \{ |z_i| < |t|^{1/2} \} \subset M_t$. We note that $M_t = M_{t,1} \cup M_{t,2}$ for each $t \in \dot{D}$ with an overlap along the vanishing $1$-cycle given by the equation $|x|=|y|=|t|^{1/2}$ in $S$. 

The next proposition deals with the limit behavior of the canonical K\"ahler form as $t \to 0$ in Fay's model.
\begin{prop} \label{musep} Let $\mu_t$ be the canonical K\"ahler form on $M_t$ and $\mu_i$ that on $M_i$ for $i=1,2$. View $M_{t,i}$ as a subset of $M_i$ for $i=1,2$. Then for $x \in M_{t,i}$, $i=1,2$, we have
\[ \mu_t(x) = \frac{h_i}{h} \mu_i(x) + o(1) \]
as $t \to 0$. Here the limit $\lim_{t\to 0} o(1)/|t|$ is equal to a $(1,1)$-form on $M_i \setminus \{p\}$ which is bounded by a constant times $|dz_i|^2|z_i|^{-2}$ on $U_i$ for $i=1,2$.
\end{prop}
\begin{proof} This is Lemma 6.9 in \cite{we}. The statement follows from the expansions in Theorem \ref{faysep} and Corollary \ref{periodsep} and the explicit formula (\ref{muv}) applied to $\mu_t$ and the $\mu_i$. 
\end{proof}
Let $k_t$ be the canonical $(1,1)$-form introduced in Section \ref{prelim} on the product $M_t \times M_t$ and let $k_i$ be that form on $M_i \times M_i$ for $i=1,2$. On a product of two manifolds we let $\pi_1,\pi_2$ denote the projections on the first and second coordinate, respectively. Continuing to view each $M_{t,i}$ as a subset of $M_i$ we let $\bar{k}_t$ be the $(1,1)$-form on $M_t \times M_t$ given by
\begin{equation} \label{defk0} \bar{k}_t = \left\{ \begin{array}{ll}
-\frac{h_2}{h} (\pi_1^* \mu_1 + \pi_2^* \mu_1) + k_1 & \textrm{on} \,\, M_{t,1} \times M_{t,1} \\
-\frac{h_1}{h} (\pi_1^* \mu_2 + \pi_2^* \mu_2) + k_2 & \textrm{on} \,\, M_{t,2} \times M_{t,2} \\
\frac{h_1}{h} \pi_1^*\mu_1 + \frac{h_2}{h} \pi_2^* \mu_2 & \textrm{on} \,\, M_{t,1} \times M_{t,2} \\
\frac{h_2}{h} \pi_1^*\mu_2 + \frac{h_1}{h} \pi_2^* \mu_1 & \textrm{on} \,\, M_{t,2} \times M_{t,1} \\
\end{array} \right. \, .
\end{equation}
We put $\bar{\nu}_t = \bar{k}_t^2$ on $M_t \times M_t$ and $\nu_i=k_i^2$ on $M_{t,i} \times M_{t,i}$ for $i=1,2$. One readily derives that
\begin{equation} \label{k0sq}  \bar{\nu}_t = \left\{ \begin{array}{ll}
-\frac{2h_2(2h-h_2)}{h^2} \pi_1^* \mu_1 \, \pi_2^* \mu_1 + \nu_1 & \textrm{on} \,\, M_{t,1} \times M_{t,1} \\
-\frac{2h_1(2h-h_1)}{h^2} \pi_1^* \mu_2 \, \pi_2^* \mu_2 + \nu_2 & \textrm{on} \,\, M_{t,2} \times M_{t,2} \\
\frac{2h_1h_2}{h^2} \pi_1^*\mu_1 \,  \pi_2^* \mu_2 & \textrm{on} \,\, M_{t,1} \times M_{t,2} \\
\frac{2h_1h_2}{h^2} \pi_1^*\mu_2 \, \pi_2^* \mu_1 & \textrm{on} \,\, M_{t,2} \times M_{t,1} \\
\end{array} \right. \, .
\end{equation}
\begin{prop} \label{ksep} Let $\nu_t=k_t^2$. We have
\[ \nu_t = \bar{\nu}_t + o(1) \]
on $M_t \times M_t$. Here $\lim_{t\to0} o(1)/|t|$ is a $(2,2)$-form on the union of the $(M_i \setminus \{p\}) \times (M_j \setminus \{p\})$ for $i,j=1,2$, bounded by a constant times $|dz_i|^2|z_i|^{-2}|dz_j|^2|z_j|^{-2}$ on $U_i \times U_j$,
by $|dz_i|^2 |z_i|^{-2}$ times a $(1,1)$-form on $M_j$ on $U_i \times (M_j \setminus U_j)$, and by a $(1,1)$-form on $M_i$ times $|dz_j|^2 |z_j|^{-2}$ on $(M_i \setminus U_i) \times U_j$, for $i,j=1,2$.
\end{prop}
As an illustration of this proposition, note that by using the identity
\[ -\frac{2h_2(2h-h_2)}{h^2} + (2-2h_1)  -\frac{2h_1(2h-h_1)}{h^2} + (2-2h_2) + \frac{4h_1h_2}{h^2} = 2 -2h \]
one sees from the expressions in (\ref{k0sq}) that the fiber integral $\int_{M_t^2} \bar{
\nu}_t$ tends to the constant fiber integral $\int_{M_t^2} 
\nu_t=2-2h$ as $t \to 0$.
\begin{proof}[Proof of Proposition \ref{ksep}] This follows from Theorem \ref{faysep}, Proposition \ref{musep} and formula (\ref{kDeltav}) applied to $k_t$ and the $k_i$. 
\end{proof}
Let $g_t$ be the canonical Green's function on $M_t$. In order to determine the asymptotic behavior of $\varphi(M_t)$ we use equation (\ref{phikDelta}) which gives that
\begin{equation} \label{phikDeltat}
\varphi(M_t) = \int_{M_t \times M_t} g_t(x,y) \, \nu_t(x,y) 
\end{equation}
for all $t$.
The discussion in Section 6 of \cite{we} leads to a suitable asymptotic expansion of $g_t(x,y)$. Following the beginning of that section we make a reparametrization of our family by putting $\tau = d_1 d_2 t$ with
\[ \log d_1 = \lim_{x \to p} \left[ g_1(x,p) - \log|z_1(x)| \right] \, , \quad
\log d_2 = \lim_{x \to p} \left[ g_2(x,p) - \log|z_2(x)| \right] \, . \]
Here $g_i$ is the canonical Green's function on $M_i$ for $i=1,2$ and $p$ still denotes the points $p_i$ on $M_i$.
\begin{prop} \label{greensep} For $t$ small enough and for local holomorphic sections $x,y$ of $M_t$ with both $x,y$ in $M_1 \setminus \{p\}$ we have the expansion
\[ g_t(x,y) = \left( \frac{h_2}{h} \right)^2 \log |\tau| + g_1(x,y) - \frac{h_2}{h}(g_1(x,p)+g_1(y,p)) + o(1) \]
as $t \to 0$. For local holomorphic sections $x,y$ of $M_t$ with $x \in M_1 \setminus \{p\}$ and $y \in M_2 \setminus \{p\}$ and $t$ small enough we have the expansion
\[ g_t(x,y) = -\frac{h_1h_2}{h^2} \log |\tau| + \frac{h_1}{h} g_1(x,p) + \frac{h_2}{h} g_2(y,p) + o(1) \]
as $t \to 0$. In both cases the limit $\lim_{t \to 0} o(1)/|t|$ is bounded in $x,y$.
\end{prop}
\begin{proof} This follows from Theorem 6.10 in \cite{we} and its proof.
\end{proof}
We are now ready to state the main result of this section.
\begin{thm} For Fay's family $M_t$ over the punctured open unit disk, the expansion
\[ \varphi(M_t) = -\frac{2h_1h_2}{h} \log|\tau| + \varphi(M_1) + \varphi(M_2) + o(1) \]
holds as $t\to 0$.
\end{thm}
\begin{proof} We study the asymptotic behavior of the fiber integral in (\ref{phikDeltat}) using the expansions of $\nu_t$ and $g_t$ in Propositions \ref{ksep} and \ref{greensep}, respectively. First of all, by dominated convergence the integral of $\nu_t$ against the $o(1)$ term in Proposition \ref{greensep} vanishes in the limit as $t \to 0$. Second, the integral of the constant term in the expansions in Proposition \ref{greensep} against the $o(1)$ term of Proposition \ref{ksep} vanishes in the limit as $t \to 0$. This is clear for both $x,y$ lying outside of the disks $U_i$. On $U_i \times U_j$, integration against the $o(1)$ term of Proposition \ref{ksep} of terms of type $\log |z_i|$ arising from the Green's functions $g_i(p,\cdot)$ yields terms of order
\[ \sim |t| \int_{|z_i|,|z_j|> |t|^{1/2}} \frac{ |dz_i|^2 |dz_j|^2}{|z_i|^2 |z_j|^2} \log|z_i| \sim (\textrm{cst.}) \cdot |t| (\log|t|)^3 \]
and these vanish indeed in the limit. A similar analysis applies to the regions of the type $U_i \times (M_j \setminus U_j)$ and $(M_i \setminus U_i) \times U_j$.

Propositions \ref{ksep} and \ref{greensep} then yield, using (\ref{normalization}), (\ref{area}), (\ref{fiberintproj}) and (\ref{k0sq}):
\begin{align*}
\int_{M_{t,1}^2} g_t \, 
\nu_t = &  \int_{M_1^2}  \left(
\left( \frac{h_2}{h} \right)^2 \log|\tau| + g_1(x,y) - \frac{h_2}{h}(g_1(x,p)+g_1(y,p)) \right) \times \\ & \times \left( -\frac{2h_2(2h-h_2)}{h^2} \pi_1^* \mu_1 \, \pi_2^* \mu_1 + \nu_1 \right) + o(1) \\
= &  \left( \frac{h_2}{h} \right)^2 \left( -\frac{2h_2(2h-h_2)}{h^2} + 2 -2h_1 \right) \log |\tau| + \varphi(M_1) + o(1)
\end{align*}
as $t \to 0$ and similarly
\begin{align*}
\int_{M_{t,2}^2} g_t \, \nu_t =  \left( \frac{h_1}{h} \right)^2 \left( -\frac{2h_1(2h-h_1)}{h^2} + 2 -2h_2 \right) \log |\tau| + \varphi(M_2) + o(1)
\end{align*}
as $t \to 0$. Next we have, using (\ref{normalization}),
\begin{align*}
\int_{M_{t,1}\times M_{t,2}} g_t \, \nu_t = &  \int_{M_1 \times M_2}   
\left( -\frac{h_1h_2}{h^2} \log|\tau| + \frac{h_1}{h}g_1(x,p)+ \frac{h_2}{h}g_2(y,p) \right) \times \\ & \times \left( \frac{2h_1h_2}{h^2} \pi_1^* \mu_1 \, \pi_2^* \mu_2  \right) + o(1) \\
= &   -\frac{h_1h_2}{h^2} \cdot \frac{2h_1h_2}{h^2} \log |\tau|  + o(1)
\end{align*}
as $t \to 0$ and by symmetry (\ref{symmetry})
\begin{align*}
\int_{M_{t,2}\times M_{t,1}} g_t \, \nu_t =   -\frac{h_1h_2}{h^2} \cdot \frac{2h_1h_2}{h^2} \log |\tau|  + o(1)
\end{align*}
as well as $t \to 0$. We obtain the required equality by adding the four integrals.
\end{proof}

\section{Non-separating case} \label{nonsep}

In this section we study Fay's model of a family of degenerating Riemann surfaces developing a single non-separating node. The basic references for this section are again \cite{fay} and \cite{we}, esp. Sections 4 and 7 of the latter. 

Let $M$ be a compact and connected Riemann surface of genus $h \geq 1$, and let $a,b$ be two distinct points in $M$. Let $z_a \colon U_a \isom D$ and $z_b \colon U_b \to D$ be local coordinate charts of $M$ centered at $a$ resp. $b$, where $D$ is the open unit disk in the complex plane. We assume that $U_a \cap U_b$ is empty. As before we let $\dot{D}$ be the punctured unit disk and we put $S = \{ (x,y,t) \in \dot{D}^3 : xy=t \}$. For each $t \in \dot{D}$ we remove from $M$ the sets given by $|z_a| \leq |t|$ and $|z_b|\leq |t|$; then we glue the annuli $\{ |t| < |z_a| <1 \}$ and $\{ |t| < |z_b| < 1\}$ on the resulting open surface along the maps given by
\[ z_a \mapsto (z_a, t/z_a,t) \in S \, , \quad z_b \mapsto (t/z_b,z_b,t) \in S \, . \]
We obtain a family of compact Riemann surfaces $M_t$ over $\dot{D}$ of genus $h+1$. This family extends to a proper analytic family $M \to D$ whose central fiber $M_0$ is the surface $M$ with the points $a,b$ identified.

We choose a symplectic basis of $\mathrm{H}_1(M_t,\zz)$ as follows. First of all we extend across $(M\setminus U_a \setminus U_b) \times \dot{D}$ a symplectic basis of homology of $M$ lying in $M\setminus U_a \setminus U_b$. We need two more loops $A_{h+1},B_{h+1}$. For the loop $A_{h+1}$ we take the boundary of the disk $U_b$. The loop $B_{h+1}$ runs across the collar connecting $U_a$ and $U_b$. Note that as $t$ traces out a generator of the fundamental group of $\dot{D}$, the collar is twisted, and the resulting loop $B'_{h+1}$ differs from the original one by $\pm A_{h+1}$. As a consequence, the period matrix $\Omega_t$ on our chosen symplectic basis is multi-valued. The next key result, cf. \cite{fay}, Proposition 3.7 describes the limit behavior of the normalized differentials of our family with the chosen homology basis.
\begin{thm} \label{faynonsep}
(Fay) For $t \neq 0$ sufficiently small, the elements of the normalized basis of holomorphic differentials $v_1(x,t),\ldots,v_{h+1}(x,t)$ of $M_t$ corresponding to the chosen homology basis have the following expansions for $i=1,\ldots, h$:
\[ v_i(x,t) = v_i(x) - t\, (v_i(b)\omega(x,a)+v_i(a)\omega(x,b)) + o(t) \, , \quad x \in M \setminus U_a \setminus U_b \, , \]
and
\[ v_{h+1}(x,t) = \frac{1}{2\pi \sqrt{-1}}\omega_{b-a}(x) -t\, (\gamma_1 \omega(x,b) + \gamma_2\omega(x,a))+o(t) \, , \quad x \in M \setminus U_a \setminus U_b \, . \]
Here the $v_i$, $i=1,\ldots,h$ are a normalized basis of holomorphic differentials on $M$, $\omega$ is the canonical differential of the second kind on $M \times M$, and $\omega_{b-a}$ is the canonical differential of the third kind on $M$ with simple poles of residues $-1,+1$ at $a,b$. The $\gamma_i$ are constants, and the evaluations at $a,b$ are carried out in the local coordinates $z_a,z_b$. Each $o(t)$ term is a holomorphic differential on $M_t$ with the property that $\lim_{t\to 0} o(t)/t^2$ is a meromorphic differential on $M$ with poles only at $a$ or $b$. 
\end{thm}
\begin{cor} \label{periodnonsep}  
The period matrix $\Omega_t$ of the $v_i(x,t)$, with $i=1,\ldots,h+1$ has the expansion
\[ \Omega_t = \left( \begin{array}{cc}
\Omega_{ij} & a_i \\
a_j & \frac{1}{2\pi \sqrt{-1}} \log t + c \end{array} \right) + o(1) \]
as $t \to 0$, where $\Omega_{ij}$ is the period matrix of $M$ with the given symplectic homology basis, $c$ is a constant, and $a_i = \int_a^b v_i$.
\end{cor}
\begin{proof} See \cite{fay}, Corollary 3.8.
\end{proof} 
Theorem \ref{faynonsep} and Corollary \ref{periodnonsep} lie at the basis of the proofs of the propositions below. Let $N_t$ be the surface $M \setminus U_a \setminus U_b$. Let $C_t$ denote the closure in $M_t$ of the `collar' $M_t \cap S$. We can view $C_t$ as the union of the annuli $\{ |t|^{1/2} \leq |z_a| \leq 1\}$ and $\{ |t|^{1/2} \leq |z_b| \leq 1\}$, with an overlap along the vanishing $1$-cycle given by the equation $|x|=|y|=|t|^{1/2}$ in $S$. We have $M_t = N_t \cup C_t$, with an overlap consisting of the two boundary cycles of $U_a$ and $U_b$. We define
\begin{align*} \mu_C = & \frac{\sqrt{-1}}{2} (\mathrm{Im} \, \Omega_t)^{-1}_{h+1,h+1} v_{h+1}(x,t) \bar{v}_{h+1}(x,t) \\
& + \frac{\sqrt{-1}}{2} \sum_{i=1}^h (\mathrm{Im} \, \Omega_t)^{-1}_{i,h+1}(v_i(x,t)\bar{v}_{h+1}(x,t) + v_{h+1}(x,t)\bar{v}_i(x,t))
\end{align*}
on $M_t$. Let $\mu_t$ be the canonical K\"ahler form on $M_t$, and $\mu$ that on $M$.  
\begin{prop} \label{munonsep}
The expansion 
\[ \mu_t =  \frac{h}{h+1} \mu  + \frac{1}{h+1} \mu_C + o(1)  \]
holds on $M_t$. Here the limit $\lim_{t \to 0}  o(1)/|t|$ is equal to a $(1,1)$-form on $M \setminus \{a,b\}$ which is bounded by a constant times $|dz_{a,b}|^2|z_{a,b}|^{-2}$ on $U_{a,b}$. 
\end{prop}
\begin{proof} This follows from the expansions in Theorem \ref{faynonsep} and Corollary \ref{periodnonsep} and the explicit formula (\ref{muv}) applied to $\mu_t$ and $\mu$.
\end{proof}
Note that the integral $\int_{N_t} \mu$ can be chosen arbitrarily close to unity by shrinking the open disks $U_a,U_b$. It follows that in the limit as $t \to 0$, the form $\mu_t$ can be written as $\frac{h}{h+1}\mu + o(1)$ on $N_t$, and as $\frac{1}{h+1}\mu_C + o(1)$ on $C_t$ (cf. \cite{jo}, proof of Lemma~4.7). This gives $N_t$ asymptotic volume equal to $\frac{h}{h+1}$ as $t \to 0$, and $C_t$ asymptotic volume equal to $\frac{1}{h+1}$.

Let $k_t$ be the canonical $(1,1)$-form introduced in Section \ref{prelim} on $M_t \times M_t$, and let $k$ be that on $M \times M$. We define
\begin{align*} d_C =& -\sqrt{-1} (\mathrm{Im} \, \Omega)^{-1}_{h+1,h+1} (v_{h+1}(x,t)\bar{v}_{h+1}(y,t) + v_{h+1}(y,t)\bar{v}_{h+1}(x,t)) \\
& - \sqrt{-1} \sum_{i=1}^h (\mathrm{Im} \, \Omega_t)^{-1}_{i,h+1}(v_i(x,t)\bar{v}_{h+1}(y,t) + v_{h+1}(x,t)\bar{v}_i(y,t)) \, . 
\end{align*}
On a product of two manifolds we denote by $\pi_1,\pi_2$ the projections on the first and second coordinate, respectively. Put
\begin{equation}
\bar{k}_t = \left\{ \begin{array}{ll}
-\frac{1}{h+1}(\pi_1^* \mu + \pi_2^* \mu) + k & \textrm{on} \,\, N_t \times N_t \\
\frac{1}{h+1} ( h \pi_1^* \mu + \pi_2^* \mu_C) & \textrm{on} \,\, N_t \times C_t \\
\frac{1}{h+1} ( \pi_1^* \mu_C + h \pi_2^* \mu) & \textrm{on} \,\, C_t \times N_t \\
\frac{1}{h+1} (\pi_1^* \mu_C + \pi_2^* \mu_C) + d_C & \textrm{on} \,\, C_t \times C_t
\end{array} \right.
\end{equation}
on $M_t \times M_t$. Let $\bar{\nu}_t=\bar{k}_t^2$ on $M_t \times M_t$ and let $\nu =k^2$ on $M \times M$. One readily derives that
\begin{equation} \label{kbarnonsep}
\bar{\nu}_t = \left\{ \begin{array}{ll}
-\frac{4h+2}{(h+1)^2} \pi_1^* \mu \, \pi_2^* \mu + 
\nu & \textrm{on} \,\, N_t \times N_t \\
\frac{2h}{(h+1)^2} \pi_1^* \mu \, \pi_2^* \mu_C & \textrm{on} \,\, N_t \times C_t \\
\frac{2h}{(h+1)^2} \pi_1^* \mu_C \, \pi_2^* \mu & \textrm{on} \,\, C_t \times N_t \\
-\frac{2h(h+2)}{(h+1)^2} \pi_1^*\mu_C \, \pi_2^* \mu_C & \textrm{on} \,\, C_t \times C_t
\end{array} \right.
\end{equation}
on $M_t \times M_t$.
\begin{prop} \label{knonsep} Let $\nu_t=k_t^2$. We have
\[ \nu_t = \bar{\nu}_t + o(1) \]
on $M_t \times M_t$. Here the limit $\lim_{t \to 0}  o(1)/|t|$ is equal to a $(2,2)$-form on $(M \setminus \{a,b\})^2$ which is bounded by a constant times $|dz_i|^2|z_i|^{-2}|dz_j|^2|z_j|^{-2}$ on $U_i \times U_j$, where $i,j \in \{a,b\}$, by $|dz_i|^2 |z_i|^{-2}$ times a $(1,1)$-form on $M$ on $U_i \times (M \setminus \{j\})$, and by a $(1,1)$-form on $M$ times $|dz_j|^2 |z_j|^{-2}$ on $(M \setminus \{i\}) \times U_j$, for $i,j \in \{a,b\}$.
\end{prop}
As an illustration of this result, using the identity
\[ -\frac{4h+2}{(h+1)^2} + 2-2h + \frac{4h}{(h+1)^2} -\frac{2h(h+2)}{(h+1)^2} = 2-2(h+1) \]
one sees from the expressions in (\ref{kbarnonsep}) that the fiber integral $\int_{M_t^2} \bar{\nu}_t$ tends to the constant fiber integral $\int_{M_t^2} \nu_t=2-2(h+1)$ as $t \to 0$.
\begin{proof}[Proof of Proposition \ref{knonsep}] This follows from Theorem \ref{faynonsep}, Proposition \ref{munonsep} and equation (\ref{kDeltav}) applied to $k_t$ and $k$. 
\end{proof}
Let $g_t$ be the canonical Green's function on $M_t$. The discussion in Section 7 of \cite{we} leads to a suitable asymptotic expansion of $g_t(x,y)$ in the case that both $x,y$ stay away from the collar. Following the beginning of Section 7 of \cite{we} we make a reparametrization of the family $M_t$ by putting $\tau = d_a d_b t$ with
\[ \log d_a = \lim_{x \to a} \left[ g(x,p) - \log|z_a(x)| \right] \, , \quad
\log d_b = \lim_{x \to b} \left[ g(x,p) - \log|z_b(x)| \right] \, . \]
Here $g$ is the canonical Green's function on $M$.
\begin{prop} \label{greennonsep}
For distinct local holomorphic sections $x,y$ of the family $M_t$ with both $x,y$ in $N_t$  the expansion
\begin{align*} g_t(x,y) = & \frac{1}{12(h+1)^2} \log|\tau| + g(x,y) + \frac{5}{6(h+1)^2} g(a,b)  \\
& - \frac{1}{2(h+1)} (g(x,a) + g(x,b) + g(y,a) + g(y,b)) + o(1)
\end{align*}
holds. Here the limit $\lim_{t \to 0} o(1)/|t|$ is bounded in $x,y$.
\end{prop}
\begin{proof} This follows from Theorem 7.2 in \cite{we}.
\end{proof}
\begin{prop} \label{fiberint} For $x$ a local holomorphic section of the family $M_t$ over a small neighborhood of zero, with $x$ lying in $N_t$ for all $t$, the expansions
\begin{align*} \int_{N_t} g_t(x,y) \mu(y) = & \frac{1}{12(h+1)^2} \log|\tau|  + \frac{5}{6(h+1)^2}g(a,b) \\ & - \frac{1}{2(h+1)} (g(x,a)+g(x,b))  + o(1)
\end{align*}
and
\begin{align*} \int_{C_t} g_t(x,y) \mu_C(y) = & -\frac{h}{12(h+1)^2} \log|\tau| - \frac{5h}{6(h+1)^2} g(a,b) \\
& + \frac{h}{2(h+1)}(g(x,a)+g(x,b)) + o(1)
\end{align*}
hold as $t \to 0$. Here $\lim_{t \to 0} o(1)/|t|$ is bounded in $x$.
\end{prop}
\begin{proof} The first equality follows from Proposition \ref{greennonsep} by dominated convergence. The second follows from the first by remarking that
\[ 0=\int_{M_t} g_t(x,y) \, \mu_t(y) = \int_{N_t} g_t(x,y) \frac{h}{h+1} \mu(y) +
\int_{C_t} g_t(x,y) \frac{1}{h+1} \mu_C(y) + o(1) \]
by (\ref{normalization}), Proposition \ref{munonsep} and the remarks following the proof of Proposition \ref{munonsep}.
\end{proof}
\begin{cor} \label{doublefiberint}
The following expansions hold:
\begin{align*}
 \int_{N_t \times N_t} g_t(x,y)\mu(x)\mu(y) & =
\frac{1}{12(h+1)^2} \log |\tau| + \frac{5}{6(h+1)^2}g(a,b) + o(1) \, , \\
 \int_{N_t \times C_t} g_t(x,y)\mu(x)\mu_C(y) & =
-\frac{h}{12(h+1)^2}\log|\tau| -\frac{5h}{6(h+1)^2}g(a,b) + o(1) \, , \\
 \int_{C_t \times C_t} g_t(x,y)\mu_C(x)\mu_C(y) & =
\frac{h^2}{12(h+1)^2} \log|\tau| +\frac{5h^2}{6(h+1)^2}g(a,b) + o(1)  \, ,
\end{align*}
as $t \to 0$.
\end{cor}
\begin{proof} The first two equalities follow from Proposition \ref{fiberint} by integrating against $\mu(x)$ over $N_t$, and dominated convergence. The third equality follows from the second by remarking that
\begin{align*} 0 = & \int_{M_t \times C_t} g_t(x,y) \mu_t(x)\frac{1}{h+1}\mu_C(y) \\  
= & \int_{N_t \times C_t} g_t(x,y) \frac{h}{h+1}\mu(x) \frac{1}{h+1}\mu_C(y) \\ & + \int_{C_t \times C_t} g_t(x,y) \frac{1}{h+1}\mu_C(x) \frac{1}{h+1}\mu_C(y) + o(1)
\end{align*}
by (\ref{normalization}), Proposition \ref{munonsep} and the remarks following the proof of Proposition \ref{munonsep}.
\end{proof}
We can now state the main result of this section.
\begin{thm} For Fay's family $M_t$ over the punctured open unit disk, the expansion
\[ \varphi(M_t) = -\frac{h}{6(h+1)} \log|\tau| + \varphi(M) - \frac{5h}{3(h+1)}g(a,b) + o(1) \]
holds as $t \to 0$.
\end{thm}
\begin{proof} Recall (\ref{phikDelta}) that $\varphi(M_t)=\int_{M_t \times M_t} g_t 
\, \nu_t$. The integral of $g_t$ against the $o(1)$-term in Proposition \ref{knonsep} vanishes in the limit as $t \to 0$. We therefore find, by Propositions \ref{knonsep} and \ref{greennonsep},
\begin{align*}
\int_{N_t \times N_t} g_t \, \nu_t = & \int_{N_t \times N_t} g_t(x,y) \cdot \left( -\frac{4h+2}{(h+1)^2} \mu(x)\mu(y) + 
\nu(x,y) \right) + o(1) \\
= &  \left( \frac{1}{12(h+1)^2} \log |\tau| + \frac{5}{6(h+1)^2}g(a,b) \right) \left(-\frac{4h+2}{(h+1)^2} + 2-2h \right) \\
& + \varphi(M) + o(1)
\end{align*}
as $t\to 0$. By Proposition \ref{knonsep} and the second equality of Corollary \ref{doublefiberint} we have
\begin{align*}
\int_{N_t \times C_t} g_t \, \nu_t = & \int_{N_t \times C_t} g_t(x,y) \frac{2h}{(h+1)^2}\mu(x)\mu_C(y) + o(1) \\
= &  \left( - \frac{h}{12(h+1)^2} \log|\tau| - \frac{5h}{6(h+1)^2} g(a,b) \right) \frac{2h}{(h+1)^2}+ o(1)
\end{align*}
as $t \to 0$ and similarly
\[ \int_{C_t \times N_t} g_t \, \nu_t =  \left( - \frac{h}{12(h+1)^2} \log|\tau| - \frac{5h}{6(h+1)^2} g(a,b) \right) \frac{2h}{(h+1)^2}+ o(1) \]
as $t \to 0$. Finally by Proposition \ref{knonsep} and the third equality of Corollary \ref{doublefiberint} we have
\begin{align*} \int_{C_t \times C_t} g_t \, \nu_t = & \int_{C_t \times C_t} g_t(x,y) \cdot - \frac{2h(h+2)}{(h+1)^2} \mu_C(x)\mu_C(y) + o(1) \\
= &  -  \left( \frac{h^2}{12(h+1)^2} \log|\tau| +\frac{5h^2}{6(h+1)^2}g(a,b) \right) \frac{2h(h+2)}{(h+1)^2}+ o(1)
\end{align*}
as $t \to 0$. Upon verifying the identities
\begin{align*} \frac{1}{12(h+1)^2} & \left(-\frac{4h+2}{(h+1)^2} +2-2h\right) - 2 \cdot\frac{h}{12(h+1)^2} \cdot \frac{2h}{(h+1)^2}  \\ &- \frac{h^2}{12(h+1)^2} \cdot \frac{2h(h+2)}{(h+1)^2} = -\frac{h}{6(h+1)}
\end{align*}
and
\begin{align*} \frac{5}{6(h+1)^2} & \left(-\frac{4h+2}{(h+1)^2} +2-2h\right) - 2 \cdot\frac{5h}{6(h+1)^2} \cdot \frac{2h}{(h+1)^2}  \\ & -\frac{5h^2}{6(h+1)^2} \cdot \frac{2h(h+2)}{(h+1)^2} = - \frac{5h}{3(h+1)}
\end{align*}
one obtains the theorem by adding all four contributions.
\end{proof}

\section{Degenerating Riemann surfaces of genus two} \label{example}

We illustrate our results in a special case. We consider families $M_t$ of Riemann surfaces of genus two, degenerating into either a union of two surfaces $M_1,M_2$ of genus one, joined at a node (separating case), or into a surface $M$ of genus one with two distinct points $a,b$ identified (non-separating case). Note that the choice of a base point on $M_1$ and $M_2$ endows both $M_1$ and $M_2$ with the structure of an elliptic curve.

We will focus on the formulas in Theorem B. We start with the separating case. Let $M_t$ be Fay's family degenerating into the union of the two elliptic curves $M_1,M_2$. We assume that $M_1,M_2$ are given as the complex tori $M_1=\cc/(\zz+\zz \omega_1)$ and $M_2=\cc/(\zz+\zz \omega_2)$, with $\omega_1,\omega_2$ elements of the complex upper half plane, and that their origins are identified. As coordinates around the origin we choose on both tori the standard euclidean coordinate $z$ coming from the uniformization by $\cc$. It follows from Corollary \ref{periodsep} that the period matrix $\Omega_t$ of $M_t$ has the expansion
\[ \Omega_t = \left( \begin{array}{cc}
\omega_1+2\pi \sqrt{-1} \, t & 2\pi \sqrt{-1} \, t \\
2\pi \sqrt{-1} \, t  & \omega_2 + 2\pi\sqrt{-1} \, t \end{array} \right) + o(t)  \]
as $t \to 0$. The first formula in Theorem B specializes  to
\begin{align} \label{firstspec} \lim_{t\to 0} & [\beta(M_t) + 4 \log|\tau|] =\\ & \nonumber \frac{2}{3}(\varphi(M_1)+\varphi(M_2)) + \frac{5}{3}(\delta(M_1)+\delta(M_2))-\frac{40}{3}\log(2\pi) \, .
\end{align}
We recall that the $\varphi$-invariant vanishes for elliptic curves. Let $\eta=\eta(\omega)$ be the Dedekind eta function,
\[ \eta = q^{1/24} \prod_{n=1}^\infty (1-q^n) \, , \quad q = \exp(2\pi \sqrt{-1}\omega) \, , \]
for $\omega$ running through the complex upper half plane, and let
\[ \|\eta\|(\omega) = (\mathrm{Im} \, \omega)^{1/4}|\eta(\omega)|   \]
be the Petersson norm of $\eta$. We note that $\|\eta\|$ is $\mathrm{SL}_2(\zz)$-invariant and hence defines an invariant of complex elliptic curves. In \cite{fa}, Section 7 it is proved that for the delta-invariant of the elliptic curve $M=\cc/(\zz+\zz \omega)$ one has
\begin{equation} \label{delta} \delta(M) = -24 \log\|\eta\|(\omega) - 8 \log(2\pi) \, .
\end{equation}
Also it follows from \cite{fa}, Section 7 that with $z$ the standard euclidean coordinate on $M=\cc/(\zz+\zz \omega)$ and for $o$ the origin on $M$, the invariant
\[ \log d = \lim_{x \to o} [g(x,o) -\log|z(x)| ] \]
evaluates to
\[ \log d = 2\log|\eta(\omega)|+\log(2\pi) \, . \]
Hence with $\tau = d_1d_2t$ and (\ref{delta}) equation (\ref{firstspec}) simplifies to
\begin{align} \label{secondspec}
\lim_{t\to 0} & [\beta(M_t) + 4 \log|t|] \\ \nonumber & = -48\log\|\eta\|(\omega_1)-48\log\|\eta\|(\omega_2) + 2\log \mathrm{Im} \, \omega_1 \, \mathrm{Im} \, \omega_2 - 48 \log(2\pi) \, .
\end{align}
We verify (\ref{secondspec}) using an explicit formula for $\beta$ in genus two. Let $\HH_2$ be the Siegel upper half space of degree two, consisting of the complex symmetric $2$-by-$2$ matrices with positive definite imaginary part. Let $\PP$ be a maximal set of pairwise inequivalent even theta characteristics in dimension two, and for each $\alpha \in \PP$ let $\theta[\alpha](z,\Omega)$ for $z \in \cc^2$ and $\Omega \in \HH_2$ be the associated theta function with characteristics,
\[ \theta[\alpha](z,\Omega) = \sum_{n \in \zz^2} \exp( \pi \sqrt{-1}  (n+a)\Omega{}^t(n+a)+2\pi \sqrt{-1} (n+a){}^t(z+b)) \, , \]
where we write $\alpha=(a,b)$ with $a,b \in (1/2)\zz^2$. We put
\[ \chi_{10}(\Omega) = \prod_{\alpha \in \PP} \theta[\alpha](0,\Omega)^2 \, , \]
the well known Siegel cusp form of weight ten on $\HH_2$. Let
\[ \|\chi_{10}\|(\Omega) = (\det \mathrm{Im} \, \Omega)^5 |\chi_{10}(\Omega)| \]
be its Petersson norm.
From \cite{djsecond}, Theorem 1.7 we obtain that
\begin{equation} \label{explicitbeta} \beta(M_t) = -2\log \|\chi_{10}\|(\Omega_t) -40\log(2\pi)+24\log 2
\end{equation}
for our family of surfaces $M_t$ with corresponding period matrices $\Omega_t$. The asymptotic behavior of $\chi_{10}(\Omega_t)$ is determined by the asymptotic behavior of the $\theta[\alpha](0,\Omega_t)$ for each of the ten $\alpha$'s. Such an analysis is carried out in, for example, \cite{dp}, Section 5.1.1. Equation (5.14) of that section yields that
\[ \chi_{10}(\Omega_t) = t^2 \, (2\pi)^4 \, 2^{12} \, \eta(\omega_1)^{24} \eta(\omega_2)^{24} + o(t^2) \]
as $t \to 0$ so that
\begin{align*}  \lim_{t \to 0} & [- \log\|\chi_{10}\|(\Omega_t) + 2 \log|t|] = \\
& -24 \log\|\eta\|(\omega_1) - 24 \log\|\eta\|(\omega_2) +\log \mathrm{Im} \, \omega_1 \, \mathrm{Im} \, \omega_2  -4\log(2\pi) - 12\log 2 \, .
\end{align*}
Combining with (\ref{explicitbeta}) we obtain an alternative derivation of (\ref{secondspec}).

Next we consider the non-separating case. Let $M=\cc/(\zz+\zz \omega)$ with $\omega$ in the complex upper half plane be an elliptic curve and choose two distinct points $a,b$ on $M$. As local coordinates around $a,b$ we choose the global euclidean coordinate $z$ given by the uniformization of $M$ by $\cc$. Let $M_t$ be Fay's family of surfaces of genus two based upon these data, degenerating into the stable curve $M$ with the points $a,b$ identified. Put $u=b-a$. It follows from Corollary \ref{periodnonsep} that the period matrix $\Omega_t$ of $M_t$ has the expansion
\begin{equation} \label{period} \Omega_t = \left( \begin{array}{cc}
\omega & u \\
u & \frac{1}{2\pi \sqrt{-1}} \log t + c
\end{array} \right) + o(1)
\end{equation}
as $t \to 0$, for some constant $c$. The second formula in Theorem B specializes to
\begin{equation} \label{firstspecnon} \lim_{t \to 0} [ \beta(M_t) + 2 \log|\tau| + 10 \log(-\log|\tau|)  ] = \frac{2}{3}\varphi(M)+ \frac{5}{3} \delta(M) - \frac{50}{3}\log(2\pi) \, .
\end{equation}
Again, the contribution $\varphi(M)$ vanishes. With (\ref{delta})
equation (\ref{firstspecnon}) becomes
\begin{equation} \label{secondspecnon}
\lim_{t \to 0}  [\beta(M_t) + 2 \log|\tau| + 10 \log(- \log |\tau|) ] = -40 \log\|\eta\|(\omega)  - 30\log(2\pi) \, .
\end{equation}
We verify (\ref{secondspecnon}) using the explicit formula (\ref{explicitbeta}) again. An analysis of the asymptotic behavior of $\chi_{10}$ on matrices of the shape in (\ref{period}) is carried out, for example, in Section 5.1.2 of \cite{dp}. Write $\omega_2$ for the $(2,2)$ entry of $\Omega_t$, and put $q = \exp(2\pi \sqrt{-1}\omega_2)$. Let $\theta=\theta(z,\omega)$ for $z \in \cc$ be the elliptic theta function with odd characteristic $(1/2,1/2)$ associated to $\omega$. Equation (5.19) of \cite{dp} yields that
\[ \chi_{10}(\Omega_t) = -q \, 2^{12} \, \eta(\omega)^{18} \, \theta(u,\omega)^2 + o(q) \]
as $t \to 0$. It follows with (\ref{explicitbeta}) that
\begin{align} \label{sub} \lim_{t\to 0} & [ \beta(M_t)  + 2 \log|q| +10\log(-\log|q|) ]  =\\
\nonumber  &- 36 \log|\eta(\omega)| - 4 \log |\theta(u,\omega)| - 10 \log \mathrm{Im} \, \omega - 30\log(2\pi) \, .
\end{align}
By Lemma 7.5 of \cite{we} we have
\[ \log |q| = \log|\tau| - 2 \, g(a,b) + 2\pi \, (\mathrm{Im} \, u)^2/\mathrm{Im} \, \omega \, , \]
where $g$ is the canonical Green's function of $M$. It is proved in Section 7 of \cite{fa} that
\[ g(a,b) = \pi (\mathrm{Im} \, u)^2/\mathrm{Im} \, \omega + \log|\theta(u,\omega)| - \log|\eta(\omega)| \, , \]
hence
\[ \log |q| = \log |\tau| -2\log|\theta(u,\omega)| + 2\log|\eta(\omega)| \, . \]
Substituting this in (\ref{sub}) we reobtain (\ref{secondspecnon}).

\vspace{0.5cm}

\noindent Address of the author:\\ \\
Mathematical Institute, \\
University of Leiden, \\
PO Box 9512, \\
2300 RA Leiden, \\
The Netherlands. \\
Email: \verb+rdejong@math.leidenuniv.nl+

\end{document}